\documentclass[11pt,reqno,oneside]{amsart}
\usepackage{verbatim,amsmath,amssymb,cite,epsfig,xspace,color,subfigure}

\numberwithin{equation}{section}

\theoremstyle{plain}

\theoremstyle{definition}

\theoremstyle{remark}

\newcommand{\ignore}[1]{}

\newcommand{\R}{\mathbb R}

\newcommand{\Her}{\mathrm{Her}}
\newcommand{\diag}{\mathrm{diag}}

\newcommand{\tr}{\operatorname{tr}}

\begin{document}

\title[Numerical solution of Dyson Brownian motion]
{Numerical solution of Dyson Brownian motion and a sampling scheme for invariant matrix ensembles} 
\author{Xingjie Helen Li and Govind Menon}





\begin{abstract}
The Dyson Brownian Motion (DBM) describes the stochastic evolution of $N$ points on the line driven by an applied potential, a Coulombic repulsion and
identical, independent Brownian forcing at each point. We use an explicit tamed Euler scheme to numerically solve the Dyson Brownian motion and sample the equilibrium measure for non-quadratic potentials. The Coulomb repulsion is too singular for the SDE to satisfy the hypotheses of rigorous convergence proofs for tamed Euler schemes~\cite{Kloeden:2012}. Nevertheless, in practice the scheme is observed to be stable for time steps of $O(1/N^2)$ and  to relax exponentially fast to the equilibrium measure with a rate constant of $O(1)$ independent of $N$. Further,  this convergence rate appears to improve with $N$ in accordance with $O(1/N)$ relaxation of local statistics of the Dyson Brownian motion. This allows us to use the Dyson Brownian motion to sample $N\times N$ Hermitian matrices from the invariant ensembles. The computational cost of generating $M$ independent samples is 
$O(MN^4)$ with a naive scheme, and $O(MN^3\log N)$ when a fast multipole method is used to evaluate the Coulomb interaction. 
\end{abstract}

\maketitle

\section{Introduction}
\subsection{The Coulomb gas on the line and Dyson Brownian motion}
Our main purpose in this article is to sample the equilibrium measure of the Coulomb gas with a given  (non-quadratic) potential $V$. The Hamiltonian for this particle system is
\begin{equation}
\label{eq:Hamiltonian}
H_N(\lambda) = \frac{1}{2}\sum_{k=1}^N V(\lambda_k) -\frac{1}{N}\sum_{1\leq j < k \leq N} \log(\lambda_j-\lambda_k).
\end{equation}
Here $\lambda =\left(\lambda_1, \ldots,\lambda_N\right)$ denotes the positions of $N$ particles on the line  and $V(s)$ is a polynomial in one variable with a term of highest degree  $\alpha_{2m}s^{2m}$ with $\alpha_{2m}>0$. We find it convenient to break permutation symmetry and work with the domain
$\mathcal{W}_N = \{\lambda \in \R^N \left| \lambda_1 < \lambda_2 < \ldots < \lambda_N \right. \}$.
At  fixed inverse temperature $\beta>0$, we define $\mu_{\beta,V}(\lambda)$ as the probability measure on $\mathcal{W}_N$ with density proportional to $e^{-\beta N   H_N(\lambda)}$. The  measure $\mu_{\beta,V}$ balances two competing effects: Coulombic repulsion pushes particles apart, but the potential $V$ confines them globally.

The Gibbs measure $\mu_{\beta,V}$ for the Coulomb gas is also the unique equilibrium measure for the Dyson Brownian motion (DBM) with potential $V$~\cite{Dyson}.  Let $B_k(t)$, $k=1, \ldots, N$ denote $N$ independent standard Brownian motions (i.e $\mathbb{E}B_k(t)=0$, $\mathbb{E}\left(B_k\left(t\right)\right)^2=t$.  The Dyson Brownian motion $\lambda(t) \in {\mathcal{W}_N}$ is the solution to the stochastic differential equation
\begin{eqnarray}
\label{eq:dyson}
\lefteqn{ d\lambda_k = - \partial_{\lambda_k} { H_N(\lambda)} \, dt + \sqrt{\frac{2}{\beta N}} \,\,dB_k,}\\
\nonumber
&&
=  \left( \frac{1}{N}\sum_{j: j \neq k} \frac{1}{\lambda_k-\lambda_j} - \frac{1}{2} V'(\lambda_k) \right)\, dt +  \sqrt{\frac{2}{\beta N}}\,\, dB_k, \quad k=1, \ldots, N.
\end{eqnarray}
This SDE is weak and strong well-posed under very general assumptions~\cite[Thm 4.3.2]{Zeitounibook}. The time scale in equation (\ref{eq:dyson}) is chosen so that a typical initial distribution relaxes to the equilibrium measure $\mu_{\beta,V}$ in time of $O(1)$. These relaxation properties underlie our simulation scheme  and allow us to rationally estimate its computational cost as follows.

We use a tamed explicit Euler scheme to numerically solve equation (\ref{eq:dyson}). (Note that a numerical solution is necessary, except when  $V(x)$ is quadratic).
Our scheme requires time steps of size $O(N^{-2})$. The naive computational cost of summing all Coulomb interactions at each time step is $O(N^2)$, but  this can be reduced to $O(N\log N)$ with a fast multipole method. Moreover, since all $\lambda_k$ are real, the domain decomposition required for the fast multipole method is very straightforward (see for example~\cite[Sec 2.1]{Greengard}). Thus, $M$ samples can be obtained with a computational effort proportional to $O(MN^3\log N)$. In the work presented here, $N$ is sufficiently small $(\leq 300)$ that we did not use a fast multipole method. The performance of our scheme (in particular, stability and speed of convergence to equilibrium) is not adversely affected by $N$. As discussed below, local convergence of DBM occurs at a rate $O(1/N)$ and this seems to {\em improve\/} the performance of the scheme as $N$ increases.  Finally, as we explain below, sampling from $\mu_{\beta,V}$ allows us to sample random Hermitian matrices from the invariant  ensembles with potential $V$.

\subsection{Invariant Hermitian ensembles}
Let $\Her(N)$ denote the space of $N \times N$ Hermitian matrices and $U(N)$ the group of $N \times N$ unitary matrices. A probability distribution $\mathbb{P}$ on $\Her(N)$ is {\em invariant\/} if for every fixed $V \in U(N)$, and every Borel set $A$ in $\Her(N)$, we have  $\mathbb{P}(M \in A) = \mathbb{P}(VMV^* \in A)$.
We focus on ensembles with a density ${ p_V(M)}$ of the form
\begin{equation}
\label{eq:invariant}
p_V(M)\, dM= \frac{1}{Z_N} e^{-N \tr V(M)} dM = \frac{1}{Z_N} e^{-N\tr V(M)} \prod_{k=1}^N dM_{kk} \prod_{j<k} dM^R_{jk} dM^I_{jk}.
\end{equation}
Here $V$ is the potential of equation~(\ref{eq:Hamiltonian}). The constant $Z_N$ is a normalizing factor
that ensures $\int_{\Her(N)}{p_V(M)} \, dM=1$, and $M_{jk}= M^R_{jk}+ \sqrt{-1}\,M^I_{jk}$ denote the real and imaginary parts of the matrix entry $M_{jk}$.

Our initial motivation for this work was to sample random matrices from the distribution $p_V(M)$ when $V$ is not quadratic (i.e. the ensemble is non-Gaussian).
Our interest in this problem stems from  observations of universal fluctuations in the time to convergence for Hamiltonian eigenvalue algorithms~\cite{PDM}. Universality was observed for a variety of Wigner ensembles, but the algorithms were not tested on matrices drawn from the invariant ensembles. The principal difficulty in sampling from an invariant ensemble is that the entries $M_{jk}$, $j \leq k$, of a random matrix $M$ are not independent (except for Gaussian ensembles). By contrast, the Wigner ensembles have independent entries and are easy to sample from. While our work arose from the desire to close this gap in~\cite{PDM}, the sampling problem is of independent interest. For example, it is well known that matrix integrals serve as generating functions for random maps~\cite{BIZ}. Thus, it is of interest to connect sampling schemes for matrices  with combinatorial algorithms for sampling random maps (e.g.~\cite{Schaeffer:2004}).

Somewhat surprisingly, despite the considerable theoretical interest in invariant ensembles~\cite{Bleher-Its,Mehta,Deift-Gioev}, there appears to have been no prior effort to numerically generate random matrices from these ensembles. In simultaneous work, Olver {\em et al\/} have also considered this  sampling problem~\cite{Trogdon}. Their approach uses numerical methods for Riemann-Hilbert problems and a constructive scheme for sampling determinantal processes, and is completely different from the methods used here. In combination, both works suggest some interesting numerical challenges arising in random matrix theory.  In particular, our work reveals the unexpected utility of Dyson Brownian motion as a test bed  for numerical methods for  stochastic differential equations.

\subsection{Determinantal structure}
\label{sec:determinantal}
The law of invariant ensembles is most easily described through the spectral representation $M = U x U^*$,  where $x=\diag(x_1,\ldots,x_N)$ denotes a diagonal matrix of (unordered) eigenvalues and $U$ is a unitary matrix of eigenvectors.  We then have
\begin{equation}
\label{eq:weyl}
e^{-N\tr V(M)} dM = e^{-N \sum_{k=1}^N V(x_k) } \triangle^2 (x)  dx \, dU.
\end{equation}
Here $\triangle(x)$ is the Vandermonde determinant $\prod_{1\leq j < k \leq N} (x_j-x_k)$, $dU$ and $dx$ denote the volume elements for normalized Haar measure on $U(N)$ and Lebesgue measure on $\mathbb{R}^N$ respectively, and $H_N(x)$ is the Hamiltonian   defined in equation (\ref{eq:Hamiltonian}).  The squared Vandermonde determinant $\triangle^2(x)$ is the Jacobian of the change of variables $M \mapsto (x,U)$~\cite{Deiftbook,Mehta}.  When included in the Hamiltonian $H_N$ as above, it has an apealing physical interpretation: the eigenvalues of invariant ensembles repel one another according to the Coulomb law and their equilibrium measure is the same as  the Coulomb gas at inverse temperature $\beta=2$. (A strict comparison requires that we reorder the $x$'s in increasing order  and include a permutation factor of $N!$ in the measure. We use $x$ to denote the unordered eigenvalues and $\lambda$ the ordered eigenvalues to avoid confusion).

The eigenvalues $x_1, \ldots, x_N$ constitute a determinantal point process on $\R$ described explicitly as follows.  Let $\pi_k(s) = s^j+\ldots$ be the monic orthogonal polynomials with respect to the measure $e^{-NV(s)}\,ds$ on $\R$, and define normalizing constants $c_k$, orthonormal functions $\phi_k$ in $L^2(\R)$, and an integral kernel $K_N$ as follows:
\begin{equation}
\label{eq:det-kernel}	
c_k^2 = \int_{\R}\pi_k^2(s)e^{-NV(s)}\, ds, \quad \phi_k(s) = \frac{1}{c_k}e^{-NV(s)/2}\pi_k(s),  \quad  K_N(r,s) = \sum_{k=0}^{N-1}\phi_k(r)\phi_k(s).
\end{equation}
Then the $m$-point correlation function  is given by
\begin{equation}
\label{eq:determinantal}
\rho_m(x_1,\ldots,x_m) = \frac{(N-m)!}{N!}\det \left(K_N(x_j,x_k)_{1 \leq j,k \leq m}\right).
\end{equation}  
(the normalization conventions in mathematics and physics appear to differ in the inclusion of the combinatorial factor of $(N-m)!/N!$ (contrast~\cite[p.108]{Deiftbook} or \cite[p.209]{Hough} with ~\cite[p.80]{Mehta})). 
The particular case $m=1$ yields the density of eigenvalues.  Under suitable assumptions on $V$ (e.g. that $V$ is convex), the limiting density of eigenvalues as $N \to \infty$ is described by a unique {\em equilibrium density\/}:
\begin{equation}
\lim_{N \to \infty} \frac{1}{N}\,  \#\{x_k \in (a,b)\} = \int_a^b { \rho(s)} \, ds, \quad -\infty <a < b < \infty.
\end{equation}
The equilibrium measure for a quartic potential $V(x)$ is described explicitly in equation (\ref{eq:quartic-eq}) below.	

\subsection{Dyson Brownian motion as a sampling method}
Our task is to sample a random matrix $M$ distributed according to equation (\ref{eq:invariant}). This is equivalent to sampling $x$ and $U$ distributed according to (\ref{eq:weyl}). Now it is immediate from (\ref{eq:weyl}) that the laws of $x$ and $U$ are independent.  Mezzadri has proposed an efficient scheme to sample uniformly distributed matrices from $U(N)$~\cite{Mezzadri}. Thus, our task is only to sample the eigenvalues according to the law (\ref{eq:determinantal}). Since the law (\ref{eq:determinantal})  is identical to the Coulomb gas at $\beta=2$, it is sufficient to solve Dyson Brownian motion until the equilibrium is approached.  This strategy is practical because of certain fundamental convergence properties of Dyson Brownian motion reviewed below.

The convergence to $\mu_{\beta,V}$ of an arbitrary initial law under DBM has been carefully studied in random matrix theory (see in particular~\cite[Sec.2]{Yau:2012} and the references therein). Assume $\lambda(t)$ solves equation (\ref{eq:dyson}), and
denote the distribution of $\lambda(t)$ by $f_t(\lambda) \mu_{\beta,V}$. Then $f_t \mu_{\beta,V}$ satifies the Fokker-Planck equation associated to equation (\ref{eq:dyson}),  and $f_t$ satisfies the evolution equation
\begin{equation}
\partial_t f_t	= \frac{1}{\beta N}\sum_{k=1}^N \partial_{\lambda_k}^2 f_t - \sum_{k=1}^N
 \partial_{\lambda_k}H (\lambda){ \partial_{\lambda_k} }  f_t.
\end{equation}
Let us assume that $V(x)$ is a convex function. Then a very minor modification of the argument in~\cite[Sec. 2]{Yau:2012} shows that there is a constant $\theta>0$ depending only on $\beta$ and $V$ such that
\begin{equation}
\label{eq:convergence}
\int_\R |f_t-1| \, d\mu_{\beta,V} \leq C e^{-\theta t},
\end{equation}
provided $f_0$ has finite relative entropy (relative to $1$) and  defines a finite Dirichlet form.

To summarize: under mild assumptions on the initial distribution, the law $f_t\mu_{\beta,V}$ converges exponentially to $\mu_{\beta,V}$ with a rate constant of $O(1)$. One of the
key recent developments in random matrix theory is a rigorous proof that the local statistics of $\lambda$ relax to equilibrium  much faster, at a time scale $O(1/N)$ (This is termed Dyson's conjecture in~\cite{Yau:2012}).
The $O(1/N)$ relaxation rate underlies the universality of fluctuations in random matrix theory. Our view is that the relaxation rates of DBM are not just of fundamental theoretical importance, but that they may be exploited in a practical sampling algorithm.  We have only used the $O(1)$ convergence in this work, but it is plausible that a more sophisticated sampling scheme could be devised using the $O(1/N)$ convergence, reducing the computational effort to $O(MN^2\log N)$.

\subsection{The numerical method}
In order to convert DBM into an effective sampling algorithm, we need a fast numerical method for the SDE (\ref{eq:dyson}). Sampling can then proceed by parallel runs with independently chosen initial data. An implicit solver for the SDE (\ref{eq:dyson}) is too expensive for this process, and it is necessary to use an explicit scheme. But here we run into an interesting complication: the well-posedness of standard explicit Euler schemes for SDE require that the vector field $\nabla_\lambda H(\lambda)$ be globally Lipschitz continuous~\cite{Higham,Kloeden-Platen}.  Moreover, counterexamples show that the Lipschitz assumption is necessary~\cite{Kloeden:2011}.
Now note that the vector field for DBM is {\em not\/} Lipschitz because of the singular nature of the Coulomb interaction.
While the counterexamples in~\cite{Kloeden:2011} do not include DBM, we have found numerically that the standard explicit Euler scheme for DBM fails to converge. This failure turns out to depend on the growth of $V'$ rather than the strong repulsion of eigenvalues and we were able to resolve it using the following tamed explicit Euler scheme that is a modification of the scheme suggested by Kloeden and his co-authors~\cite{Kloeden:2012}.

Let $\triangle t$ denote the time step, $T>0$ a fixed time of $O(1)$,  and $\lambda^n =(\lambda_1^n,\lambda_n^2, \ldots, \lambda_N^n)$ the discretized approximation to the solution $\lambda(t)$ to (\ref{eq:dyson}) at the time  $t= n \triangle t$. Assume also that $\triangle B_k^n$ are independent standard (mean zero, variance one) normal random variables for $k=1, \ldots, N$, $n=0, \ldots, T/\triangle t$. We use the following numerical scheme:
\begin{eqnarray}
	\label{eq:scheme}
	\lambda_k^{n+1}&=&\lambda_k^n +
\left(\sum_{j\neq k}\frac{1}{\lambda_n^k-\lambda^n_j}- \frac{V'(\lambda_k^n)}{2+\triangle t \left|V'(\lambda_k^n)\right|} \right) { \triangle t}
 +	\sqrt{\frac{2\triangle t }{\beta N}}\,\,   \triangle B_k^n, \\
	   \lambda_k^0 &= & \xi_k.
\end{eqnarray}
The index $k$ denotes a  coordinate in $\mathcal{W}_N$ and runs from $1$ to $N$. The index $n$ is the time step and runs from $0$ to $T/\triangle t$.  The initial data $\lambda^0_k=\xi_k$ are assumed given for the purpose of solving the SDE. In practice, they will be chosen at random as described below.

The main points to note are: (i) this scheme is explicit; (ii) the vector field $V'(x)$ has been ``tamed''; and (iii) the Coulomb interaction is left unchanged. Note also that
\begin{equation}
	\label{eq:scheme2}
	\lambda_k^{n+1}=\lambda_k^n + \left( \sum_{j\neq k}\frac{1}{\lambda_n^k-\lambda^n_j} -  \frac{1}{2} V'(\lambda_k^n)    \right) \triangle t +  \sqrt{\frac{2\triangle t }{\beta N}}\,\,   \triangle B_k^n  +    \frac{V'(\lambda_k^n) \left|V'(\lambda_k^n)\right|}{4+2\triangle t \left|V'(\lambda_k^n)\right|}\,(\triangle t)^2.
\end{equation}
Therefore, this method agrees with the explicit Euler scheme up to second-order. Since the explicit Euler scheme is consistent with the SDE, so is the tamed scheme. Thus, stability of the scheme is enough to establish convergence. However, as a consequence of the Coulomb interaction, it does {\em not\/} satisfy the assumptions of the main theorem in~\cite{Kloeden:2012}. It is not clear to us at present, if the techniques used to establish well-posedness of the SDE (\ref{eq:dyson}) can also be used to establish stability of the numerical scheme.

\section{Numerical Experiments}\label{section:experiments}
\subsection{Equilibrium measures for quartic ensembles}  
In this section we present the result of numerical experiments with the scheme (\ref{eq:scheme}) for the quartic potentials
\begin{equation}
	\label{eq:potential}
	V(x) = \frac{q x^2}{2} +  \frac{g x^4}{4}, \quad q, g \geq 0. 
\end{equation} 
The quartic ensembles are the simplest non-Gaussian invariant ensembles and provide important insights into general phenomena for invariant ensembles. For example, quartic ensembles were used to enumerate random quadrangulations~\cite{BIZ}, and to establish the first universality theorems for non-Gaussian ensembles~\cite{Bleher-Its}. Thus, we expect  that they have similar utility for a numerical study. 

In all that follows, we assume $\beta=2$ so that the samples $\lambda$ correspond to the eigenvalues of a Hermitian matrix.   If $g=0$ we obtain a Gaussian ensemble. For $q>0$ we may always rescale so that $q=1$. When $q=1$ and $0  \leq  g < \infty$, the equilibrium density for the quartic ensemble is given explicitly by~\cite{BIZ}
\begin{equation}
	\label{eq:quartic-eq}
	\rho(s) = \frac{1}{2\pi}\left(1 + 2ga^2 + gs^2\right)\sqrt{4a^2 -s^2} \, \mathbf{1}_{|s|\leq 2a}, \quad a= \sqrt{\frac{\sqrt{1+12g}-1}{6g}}.
\end{equation}
In the limit $g=0$ this density reduces to Wigner's semicircle law.
When $q=0$ the potential $V(x)$ is purely quartic and has equilibrium density
\begin{equation}
	\label{eq:quartic-pure}
	                        \rho(s)= \frac{1}{2\pi} \left(2g a^2 + gs^2\right)\sqrt{4a^2 -s^2} \, \mathbf{1}_{|s|\leq 2a}, \quad a=\left(3g\right)^{-1/4}.
\end{equation}

\subsection{Description of numerical experiments}
We fix a potential (i.e. the parameters $g$ and $q$ in equation (\ref{eq:potential})  and perform $M$ independent trials. Each trial takes the following form:
\begin{enumerate} 	
\item Choose an initial vector $\lambda^0$ at random independent of other trials.
\item Solve (\ref{eq:scheme}) for $\lambda^n$ for each $n=1,\ldots, T/\triangle t$.
\end{enumerate}
The initial vector $\lambda^0$ is typically chosen to be the vector of eigenvalues of a Gaussian Hermitian matrix. This choice is not entirely necessary, since
the universality results show that the spectrum of a random matrix from a Wigner ensemble that satisfies the fourth-moment condition $\mathbb{E}(|M_{jk}|^4) \leq C < \infty$ will relax in time of $O(1)$ to the equilibrium density $\mu_{\beta,V}$ for DBM~\cite{Yau:2012}. However, we find that there is a distinct difference in the time to equilibration between a vector $\lambda^0$ that is obtained as the eigenvalues of a random matrix, or a vector $\lambda^0$ consisting of independent entries $\lambda^0_k$, $k=1,\ldots,N$.

The typical end time is taken to be $T=O(1)$ and all simulations are observed to equilibrate when the time step is chosen to be $\triangle t = 1/N^2$. Thus, each run involves $N^2$ evaluations of the vector field, and each such evaluation costs $O(N^2)$, with a total computational cost of $O(N^4)$. As we have remarked above, this can be reduces to $O(N^3\log N)$ with a fast multipole method.  There is a further computational cost associated to the initialization, since we must compute the eigenvalues of the initial matrix. However, this cost is of lower order if the initial matrix is chosen from a tridiagonal ensemble whose eigenvalues have the law $e^{-\beta\tilde{H}_N(x)}$ where $\tilde{H}_N$ is the Hamiltonian corresponding to $V(x)=qx^2$ (for example, the $\beta$-ensembles of Dumitriu and Edelman~\cite{Dumitriu}). The total cost for $M$ trials is $O(MN^4)$.

For $M$ parallel trials, we denote by $\lambda^{n,j}_k$ the $k$-the coordinate of the solution to equation (\ref{eq:scheme}) at time step $n$ for the $j$-th trial. The  empirical distribution function  is defined to be
\begin{equation}
\label{eq:empirical-defn}
F_{n,M}(\lambda) = \frac{1}{NM} \sum_{j=1}^M \sum_{k=1}^N \mathbf{1}_{\lambda^{n,j}_k \leq \lambda},
\end{equation} 
and the exact distribution function and equilibrium distribution function are denoted
\begin{equation}
\label{eq:F-exact}
F_N(\lambda) = \frac{1}{N}\int_{-\infty}^\lambda K_N(s,s)\, ds, \quad F_\infty(\lambda) =  \int_{-\infty}^\lambda \rho(s)\, ds.
\end{equation}
We use the Kolmogorov-Smirnov (KS) statistic  
\begin{equation}
	\label{eq:KS-distance}
	D_{n,M} = \sup_\lambda \left| F_{n,M}(\lambda) - F_N(\lambda)\right|
\end{equation}	
to measure the distance between the empirical measure $F_{n,M}$ and the exact distribution $F_N$. For large enough $N$ (in practice, $N \geq 30$) we replace $F_N$ with $F_\infty$ in the formula since the error $\|F_N-F_\infty\|_{\infty}$ is much smaller than $D_{n,M}$. This is discussed below.
    
The above comparison only contrasts $1$-point statistics of the empirical and exact distributions of $x(t)$. In order to study convergence of a statistic that truly reflects the determinantal nature of the process, we compare the exact probability  $A(\theta)$ that there are no points $x_k(t)$ in the interval $(-\theta,\theta)$ with the empirical measure of $A(\theta)$  (our notation follows~\cite[Sec. 5.4]{Deiftbook}). The gap probability $A(\theta)$ is given by the Fredholm determinant \begin{equation}
\label{eq:An}
A(\theta) = \det(I-K_N),
\end{equation} 
with kernel $K_N$ acting on $L^2(-\theta,\theta)$. We compute this Fredholm determinant using a  numerical procedure introduced by Bornemann~\cite{Bornemann:2009}.

\begin{figure}[h]
\begin{center}
\subfigure[$V(x)=x^2/2+x^4/4$]{
\includegraphics [height=5 cm, width=5 cm]
{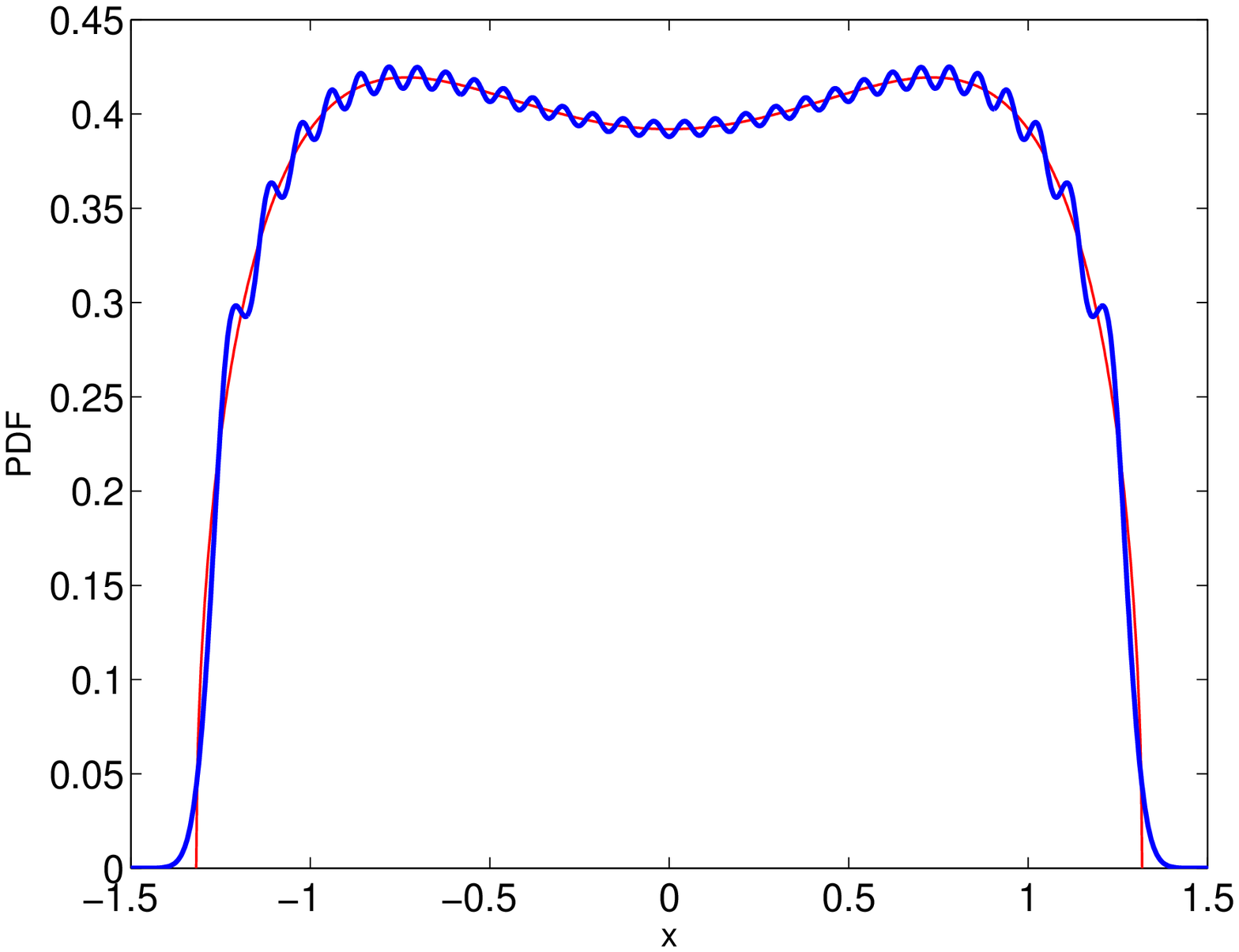}
}
\subfigure[$V(x)=x^4/4$]{
\includegraphics [height=5 cm, width=5.2 cm]
{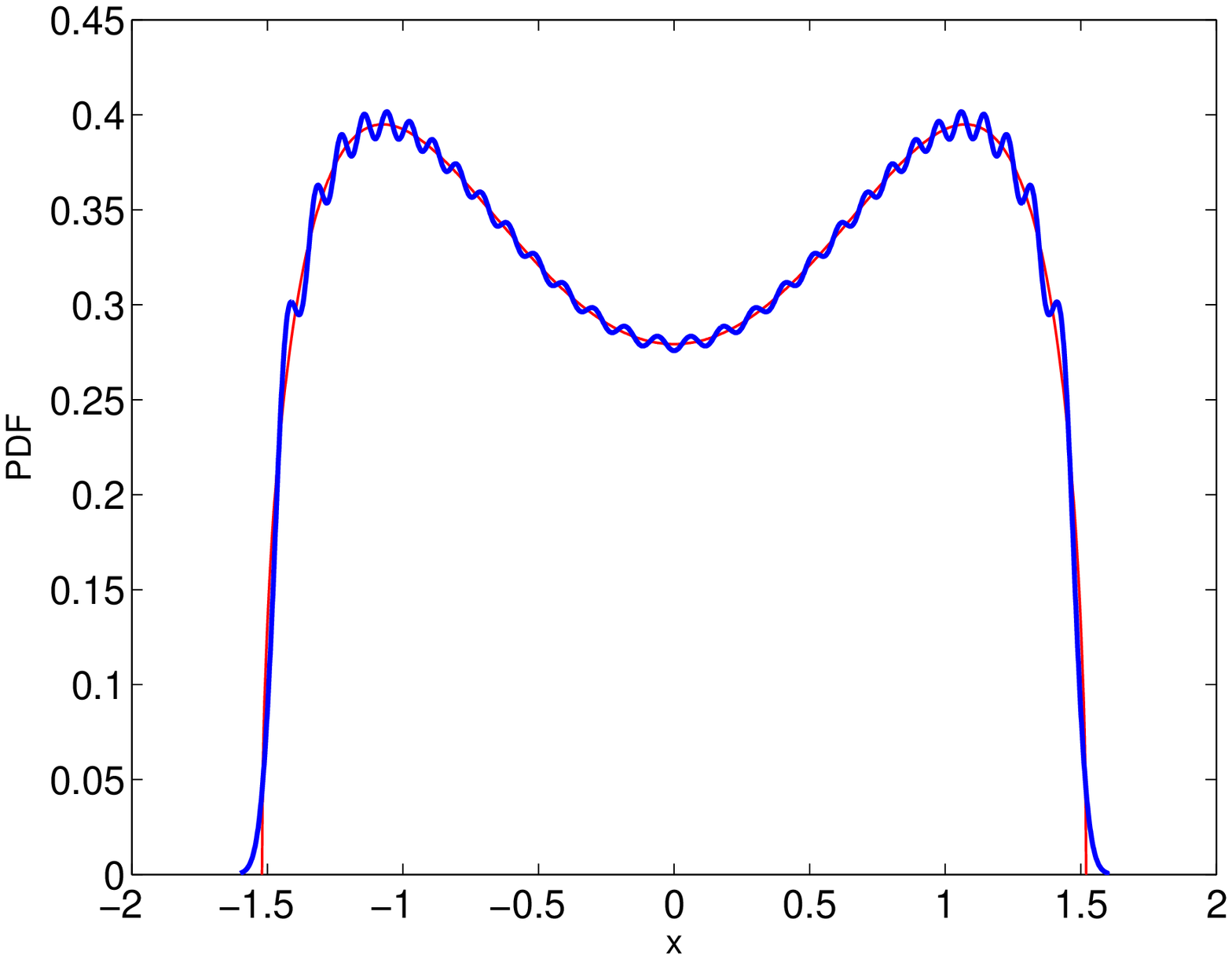}
}
\subfigure[$V(x)=x^4/4$]{
\includegraphics [height=5 cm, width=5.2 cm]{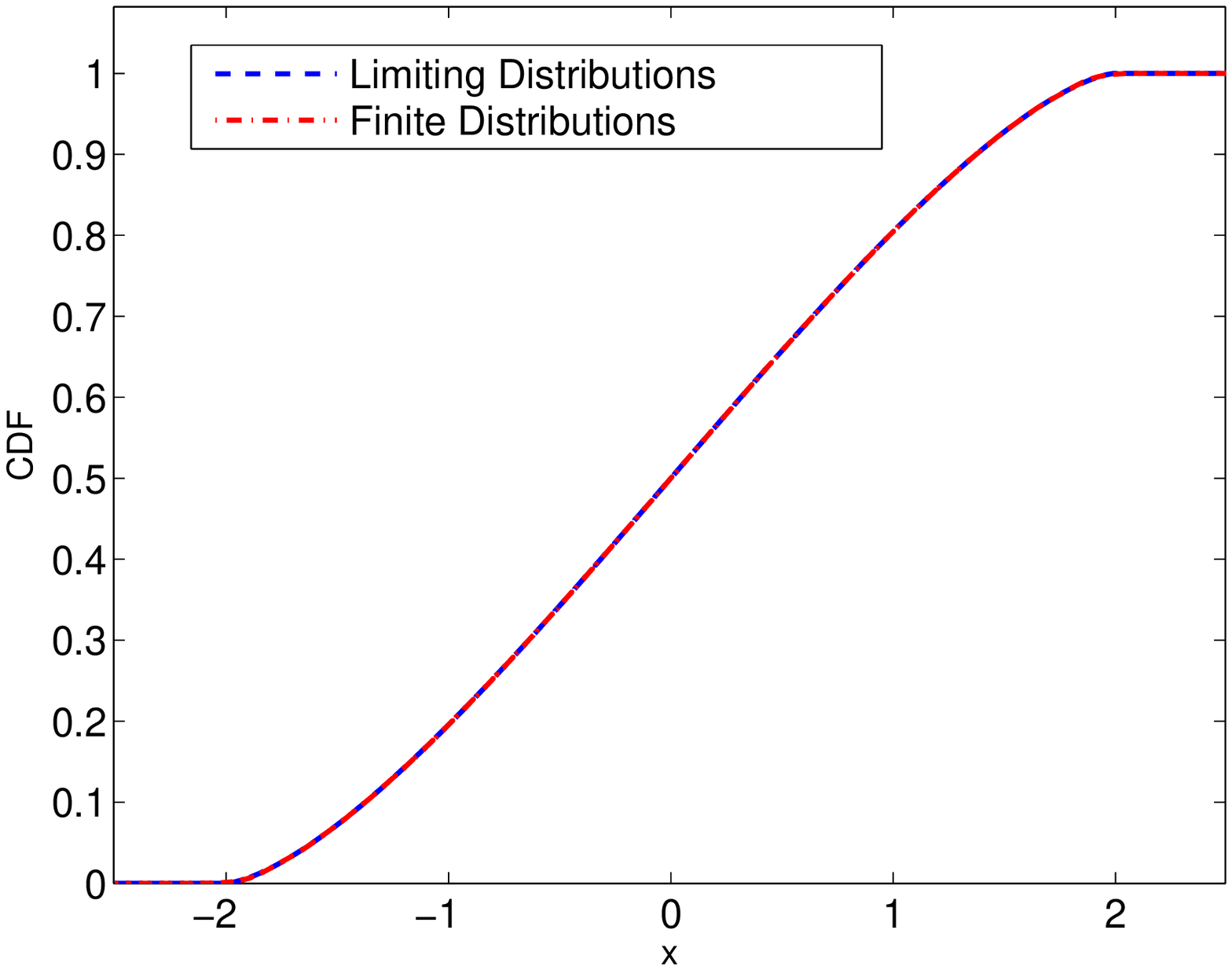}
}

\end{center}
\caption{{\bf Comparison of finite laws for $N=30$ with limiting laws}.  {\bf (a), (b)} Comparison of probability density functions.  
{\bf (c)} Cumulative-distribution functions.}\label{Fig:finitePDF}
\end{figure}

There are some numerical subtleties in these calculations.
While equation (\ref{eq:determinantal}) provides a complete prescription for the exact distribution of all correlation functions once $V$ is given, in order to  compute $F_N$ we must determine $\pi_k$, $c_k$ and $K_N$ accurately. This is a delicate problem, since $c_k$ grows rapidly. For the choices of $g$ and $q$ tested here, we used the orthogonal polynomial generator developed by Gautschi along with a standardized C++ implementation developed by Yang {\em et al}~\cite{Gautschi:1996,Yang:2012}. We find that $K_N$ can be computed accurately only for $N \leq 30$. However, as $N$ increases,  both the density $K_N(\lambda,\lambda)$ and distribution function $F_N(\lambda)$ converge exponentially fast to the density $\rho(\lambda)$ and distribution $F_\infty(\lambda)$ respectively, and the convergence is uniform in $\lambda$. The convergence of densities $K_N(\lambda,\lambda)$ to $\rho$ is illustrated  in Figure~\ref{Fig:finitePDF} for different potentials. For example, the $L^{\infty}$ distance between the finite and limiting distributions shown for $N=30$ is only $1.024e-3$.  For the same reason, we restrict our comparison of $A(\theta)$ with the empirical measure  to $N \leq 30$ (Note that our sampling scheme is not restricted to $N\leq 30$: the point is that it is numerically subtle to  compute {\em exact} formulas in random matrix theory! See~\cite{Bornemann:2009,Bornemann:2010} for a lively discussion of this issue).

\subsection{Results}
\begin{figure}[h]
\begin{center}
\subfigure[$V(x)= x^2/2+ x^4/4$]{\includegraphics [ height =7.5 cm, width=8
cm]{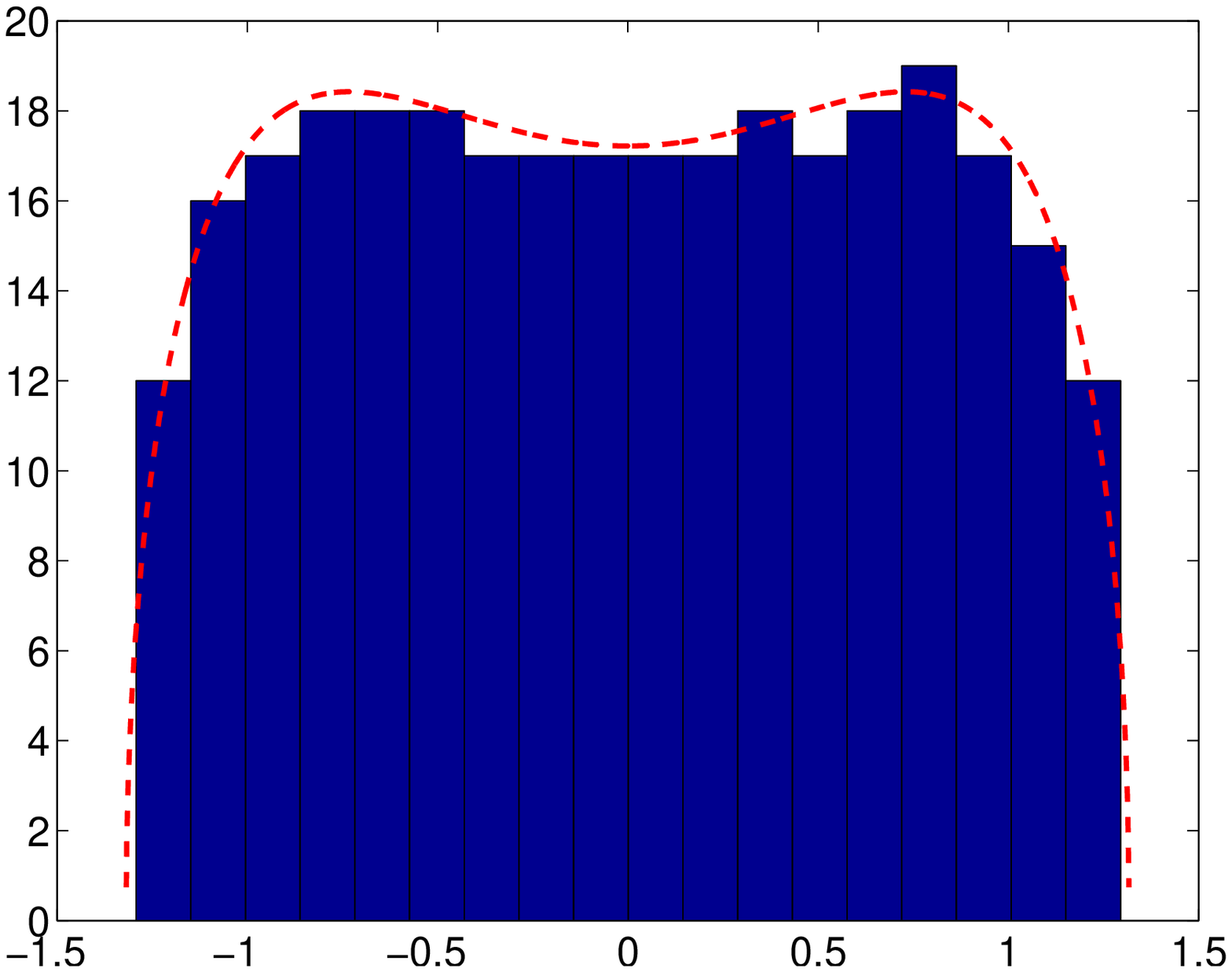}}
\subfigure[$V(x)= x^2/2+ 10 x^4/4$]{\includegraphics [ height =7.5 cm, width=8
cm]{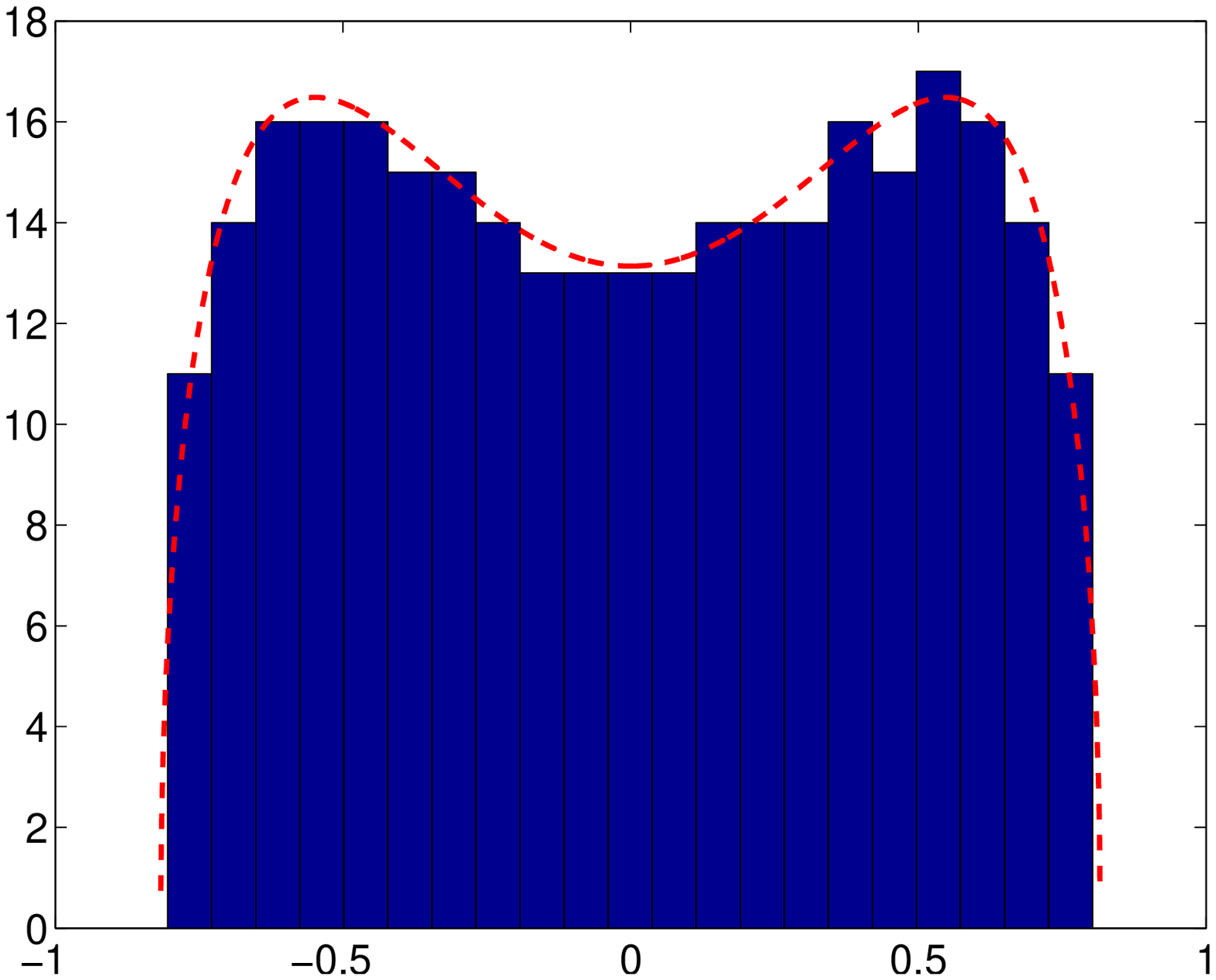}}
\end{center}
\caption{{\bf Histogram of empirical distribution after one long run.} 
A histogram of the empirical distribution of $\lambda_k$ is compared  with the theoretical limiting law (dashed line). The system size is $N=300$, the maximal time is  $T=1200$, the time step is $\triangle t=1/N^2$. }\label{Fig:hqquartic}
\end{figure}

We first present a visual comparison of the empirical distributions for a single trial with a large matrix ($N=300$). In these experiments, the time $T$ is $O(N)$. This is a much longer time scale than required, but it serves to establish that the scheme is stable. The limiting spectral distributions $F_\infty$ is compared with the empirical distribution for two different potentials in Figure~\ref{Fig:hqquartic}.  


\begin{figure}[h]
\begin{center}
\subfigure[$V(x)=x^2/2+x^4/4$ ]{\includegraphics [  height=7.5 cm, width=8
cm]{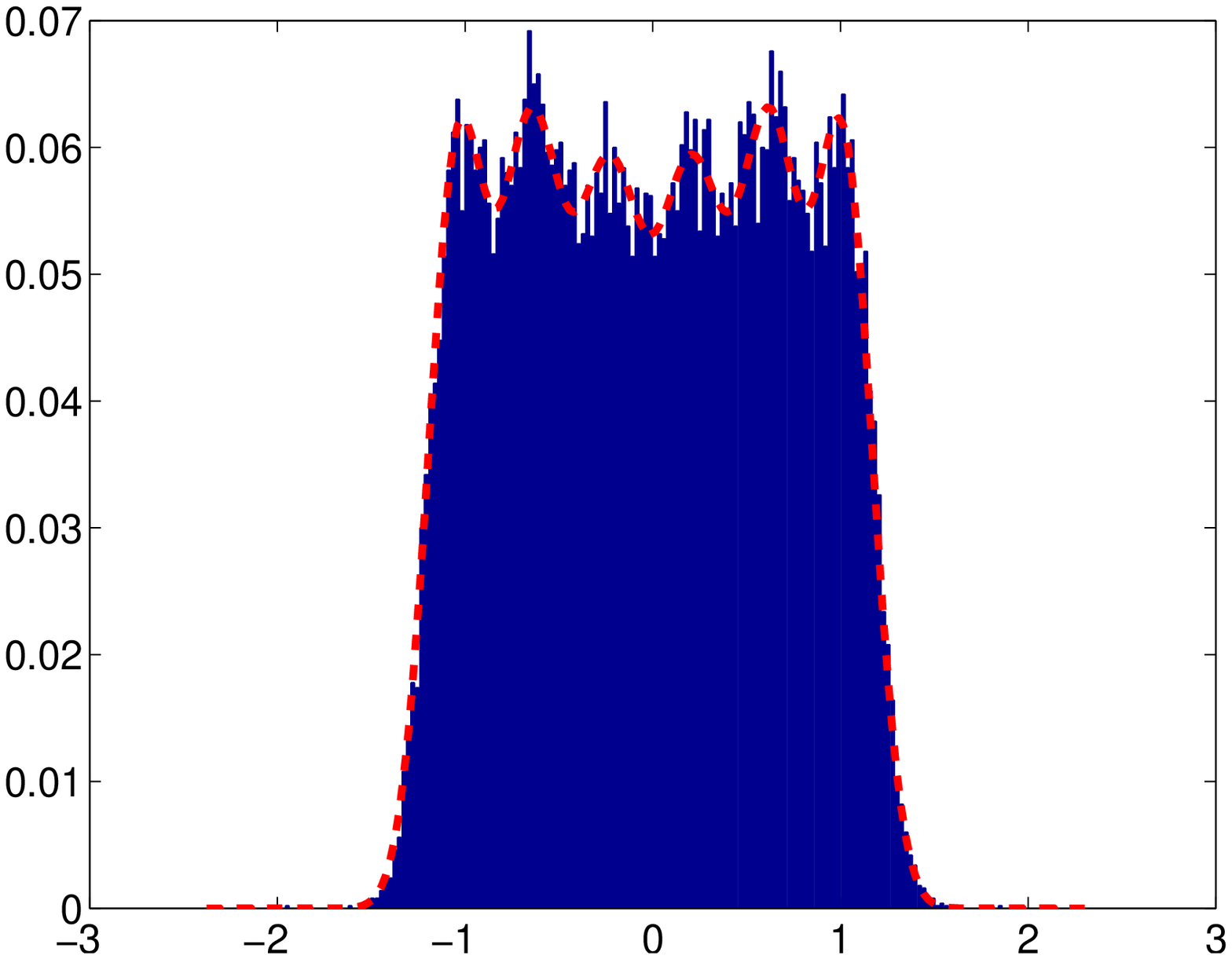}}
\subfigure[$V(x)=x^4/4$ ]{\includegraphics [  height= 7.5 cm,  width=8
cm]{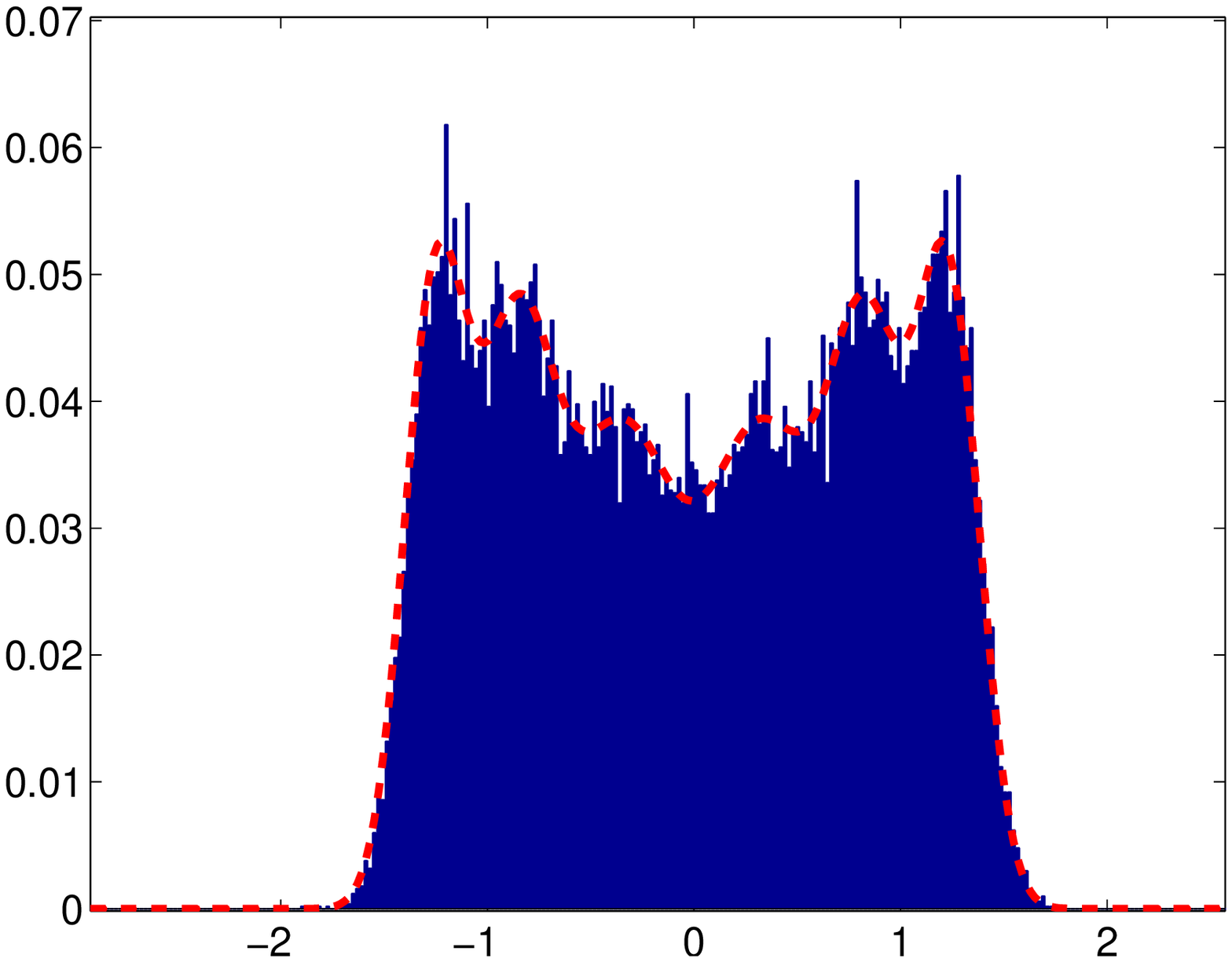}}
\end{center}
\caption{{\bf Oscillations for small $N$}. A histograms of the empirical distribution of $\lambda_k$ is compared with the exact law for two different quartic  potentials. The number of particles is $N=6$, the maximal time is $T=24$, the time step is $\Delta t=1/N^4$ and $M=5000$ trials are performed.}\label{Fig:nongoefiniteLaw}
\end{figure}

\begin{figure}[h]
\begin{center}
\includegraphics[height=8 cm, width = 11 cm]{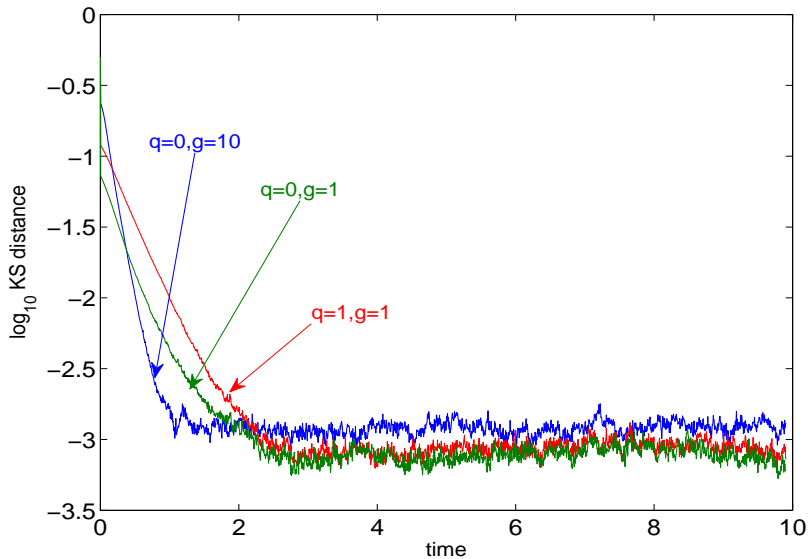}
\end{center}
\caption{{\bf Exponential decay of the KS distance}. 
The KS distance defined in equation (\ref{eq:KS-distance})  is plotted as a function of time for several ensembles. In these simulations, $N=100$, the maximal time is $T=10$, the time step is $\triangle t=1/2N^2$, and $M=1000$ trials were performed.  Exponential decay to the equilibrium distribution is seen in a time of $O(1)$.} \label{Fig:FinitevsLimit}
\end{figure}

\begin{figure}[h]
\begin{center}
\subfigure[Variation with $M$]{\includegraphics [  height=7.5 cm, width=8
cm]{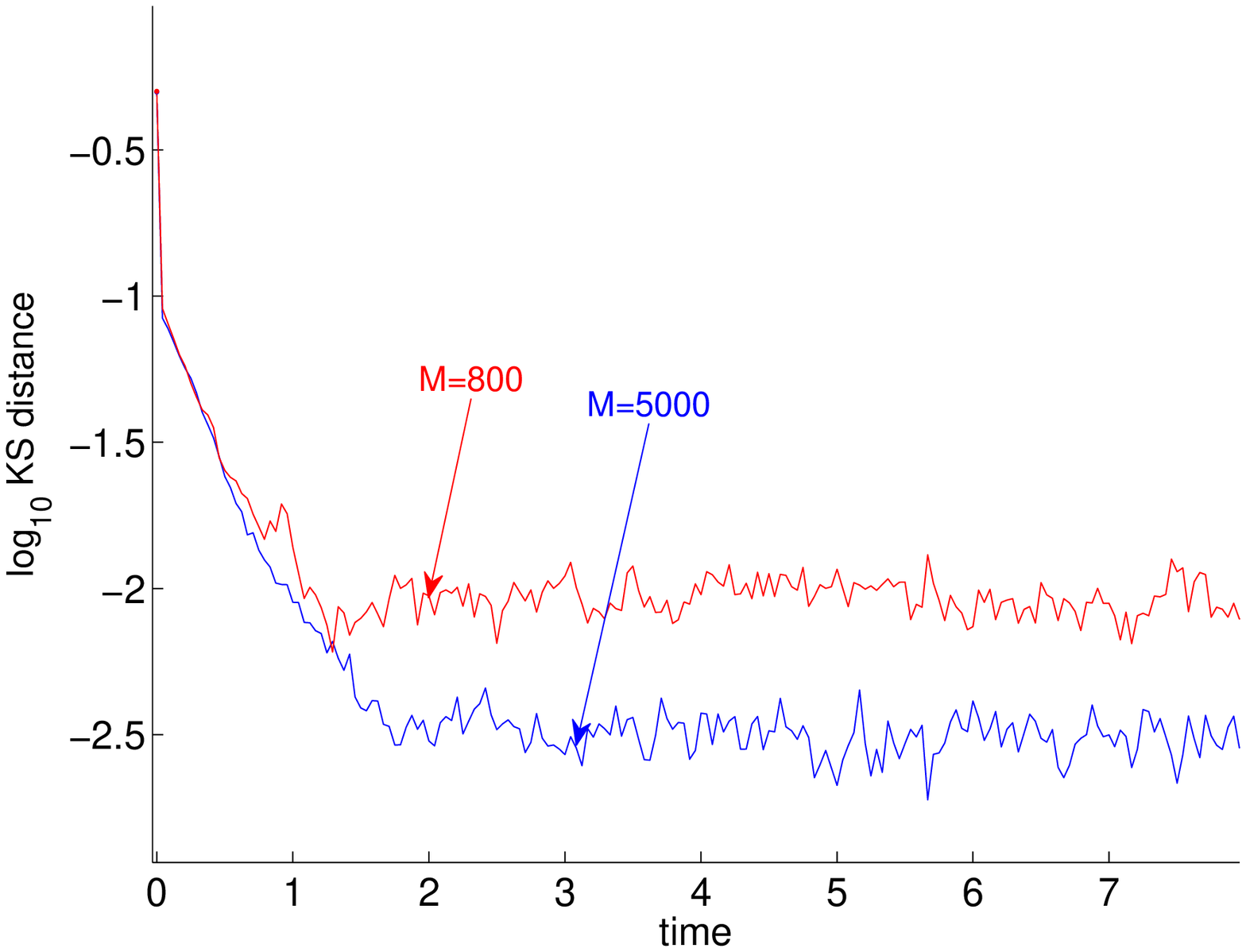}}
\subfigure[Variation with $N$ ]{\includegraphics [  height= 7.5 cm,  width=8
cm]{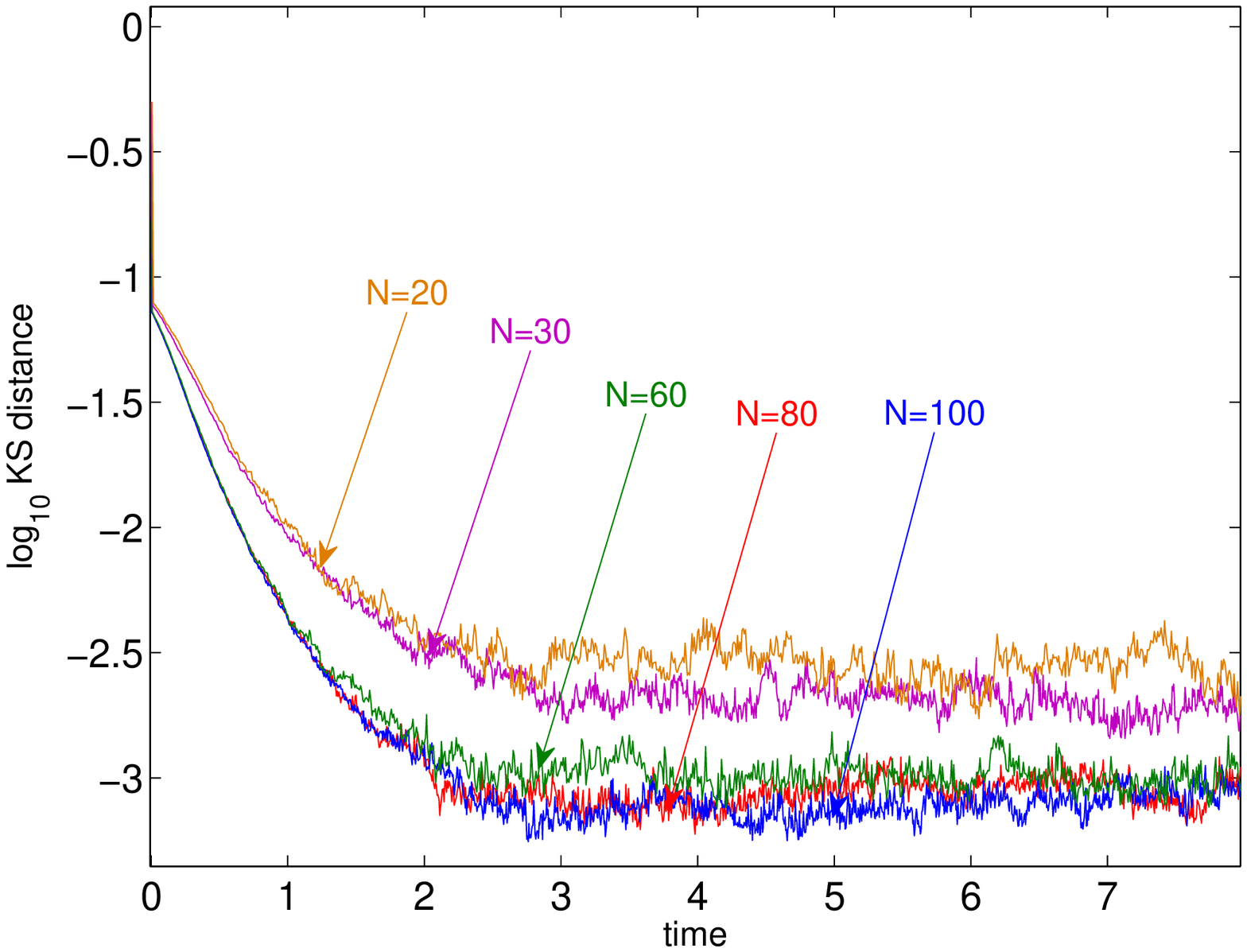}}
\end{center}
\caption{{\bf Exponential decay of the KS distance for the empirical measure}. 
The KS distance defined in equation (\ref{eq:KS-distance})  is plotted as a function of time for the ensemble with $V(x)=x^4/4$ for different values of $N$.
{\bf (a)} $N=6$, maximal time $T=8$ and time step $\triangle t=1/N^4$.
{\bf (b)} Maximal time $T=8$, time step $\triangle t=1/4N^2$ and $M=1000$ trials. }
\label{Fig:quarticfiniteErr}
\end{figure}

\begin{figure}[h]
\begin{center}
\subfigure[$A(\theta)$]{\includegraphics [  height=7.5 cm, width=8
cm]{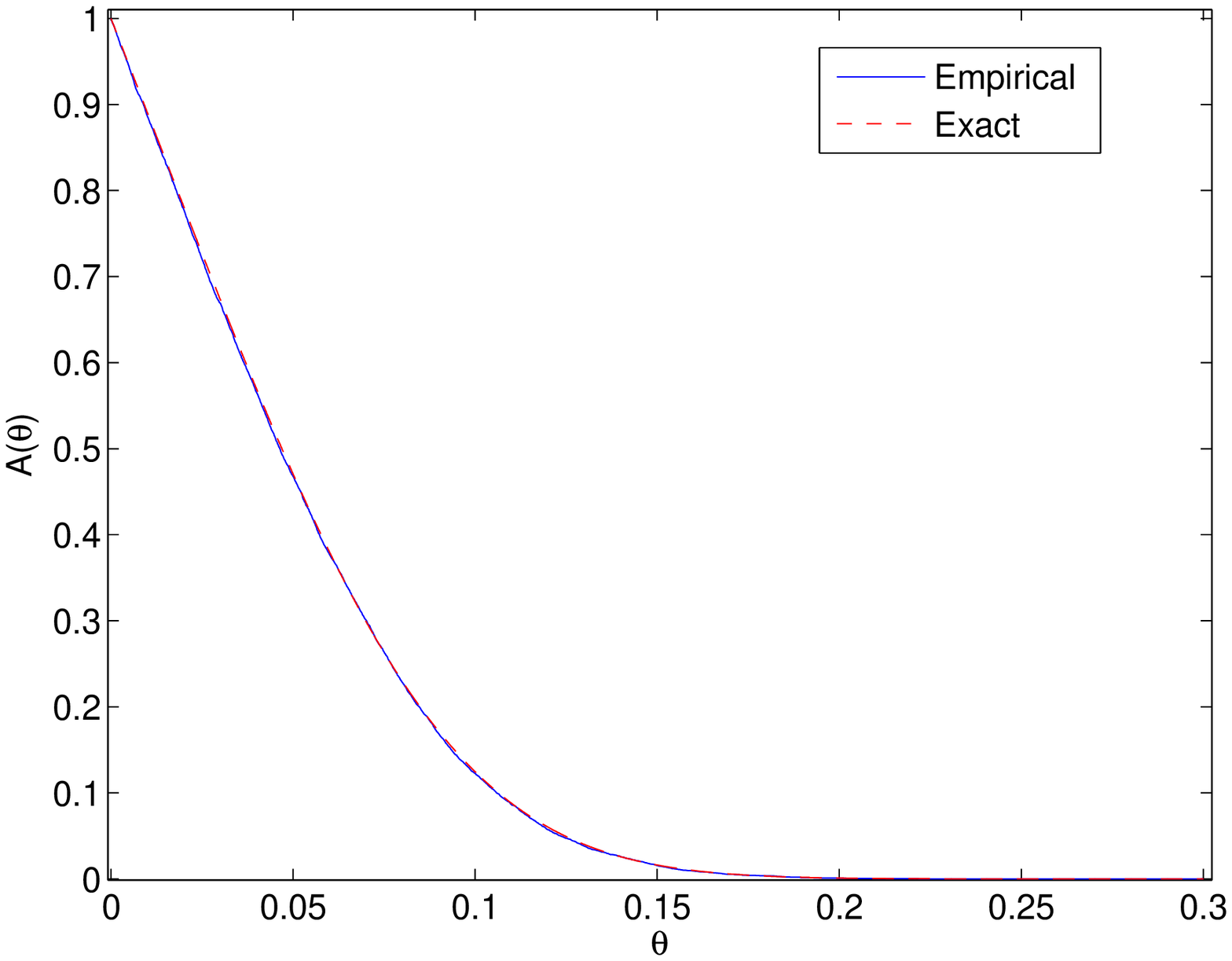}}
\subfigure[Decay of KS distance]{\includegraphics [  height= 7.5 cm,  width=8
cm]{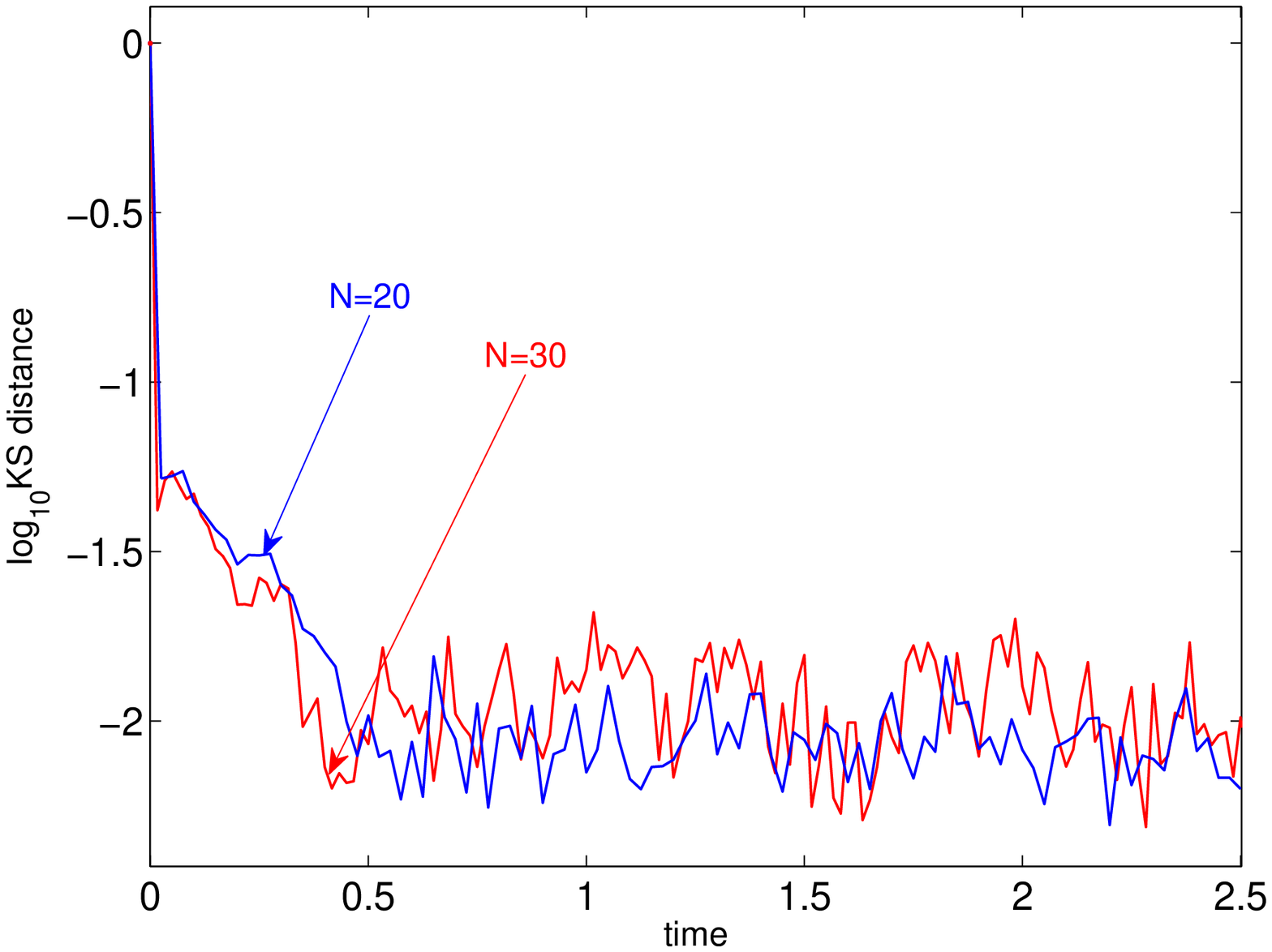}}
\end{center}
\caption{{\bf Exponential decay of the KS distance for gap probabilities.}
$A(\theta)$ for the quartic potential $V(x)=x^4/4$.
{\bf (a)} System size $N=20$, maximal time $T=2.5$ and time step $\triangle  t=1/4N^2$.
{\bf(b)} Maximal time $T=2.5$, time step
$\triangle t=1/4N^2$ and  $M=10000$ trials.}\label{Fig:quartic_AthetaErr}
\end{figure}

For small $N$ the density $\rho_N(s)$ oscillates and it is interesting to test if the scheme is sensitive to these oscillations. Of course, for small $N$ it is necessary to use a large number of trials to obtain a large enough data set.  Figure~\ref{Fig:nongoefiniteLaw}
presents the results of such a numerical experiment. We see that the scheme captures the oscillations of the density quite well. We measured the rate of convergence to equilibrium in the KS distance (\ref{eq:KS-distance}) for many different ensembles (i.e. potentials $V$), systems sizes $N$, number of samples $M$, and time steps $\triangle t$. A representative sample of the results is provided in Figure~\ref{Fig:FinitevsLimit}, Figure~\ref{Fig:quarticfiniteErr}, and  Figure~\ref{Fig:quarticfiniteErr}. We draw the reader's attention to the robust exponential decay to the equilibrium measure with a rate of $O(1)$ in all these figures. Exponential convergence to equilibrium with a rate of $O(1)$ is also seen for the gap probability $A(\theta)$ in Figure~\ref{Fig:quartic_AthetaErr}. 


\section{Conclusion}
We have used a tamed explicit Euler scheme for Dyson Brownian motion to sample the equilibrium  measure of the Coulomb gas on the line with a non-Gaussian potential $V$.
Despite the singular nature of the Coulomb interaction and the nonlinear growth of the confining potential $V$, our scheme was always observed to be stable and convergent  for time steps of $O(1/N^2)$, and to reach equilibrium in a time of $O(1)$.     
The simulation results agree with the exact laws for both finite and limiting cases. The KS distances between the empirical and exact distributions decay exponentially with respect to simulation time. We view this as a demonstration of the  practical utility of tamed explicit Euler schemes in situations where theoretical convergence results are not known to apply.

Our method can be combined with sampling techniques for matrices from the unitary group to sample a random matrix from the invariant ensemble with potential $V$. We find that the cost of $M$ trials is $O(M N^3)$, though further acceleration may be provided by the fast local relaxation of the Dyson Brownian motion. To the best of our knowledge, our work and the simultaneous work of Olver {\em et al}~\cite{Trogdon},  present the first systematic sampling schemes for invariant ensembles.

\section{Acknowledgments}
We thank Percy Deift for several stimulating discussions on numerical computations in random matrix theory. We thank Xiu Yang for assistance with the methods of~\cite{Yang:2012}.  This work has been supported in part by NSF grant DMS 07-48482.

\bibliographystyle{abbrv}
\bibliography{invsamp}
\end{document}